\documentclass{amsart}
\usepackage[T1]{fontenc}
\usepackage{amssymb}
\usepackage[hyperref]{degt}

\def\GC{\Cal{C}}
\def\num#1{\bold#1}
\let\fixed=Q
\let\octad=O
\def\set{\bar\octad}
\def\dual{^\vee}
\let\tV=V
\def\bV{\bar V}
\let\tX=X
\def\bX{X_4}
\def\Nmax#1{N_{#1}}

\let\3\relax

\title{800 conics in a smooth quartic surface}

\author{Alex Degtyarev}

\address{%
Department of Mathematics\\
Bilkent University\\
06800 Ankara, TURKEY}

\email{
degt@fen.bilkent.edu.tr}

\thanks{%
The author was partially supported by the T\"{U}B\DOTaccent{I}TAK grant 118F413%
}



\begin{document}

\begin{abstract}
We construct an example of a smooth spatial quartic surface that contains 800
irreducible conics.
\end{abstract}

\maketitle

\section{Introduction}\label{S.intro}

This short note was motivated by Barth, Bauer~\cite{Barth.Bauer:conics},
Bauer~\cite{Bauer:conics}, and my recent paper~\cite{degt:conics}.
Generalizing~\cite{Bauer:conics}, define $\Nmax{2n}(d)$ as the maximal number of
smooth rational curves of degree~$d$
that can lie in a smooth degree~$2n$ $K3$-surface
$X\subset\Cp{n+1}$. (All algebraic varieties considered in this note are over~$\C$.)
The bounds $\Nmax{2n}(1)$ have a long history and currently are well
known, whereas for $d=2$ the only known value is $\Nmax6(2)=285$
(see~\cite{degt:conics}). In the most classical case $2n=4$ (spatial quartics),
the best known examples have $352$
or $432$ conics (see~\cite{Barth.Bauer:conics,Bauer:conics}),
whereas the best known upper
bound is $5016$ (see~\cite{Bauer:conics}, with a reference to S.~A.~Str{\o}mme).

For $d=1$, the extremal configurations (for various values of~$n$) tend to exhibit
similar behaviour. Hence, contemplating the findings of~\cite{degt:conics},
one may speculate that
\roster*
\item
it is easier to count \emph{all} conics, both irreducible and reducible,
{\3but}
\item
nevertheless, in extremal configurations all conics are irreducible.
\endroster
On the other hand, famous \emph{Schur's quartic} (the one on which
the maximum $\Nmax4(1)$ is
attained) has $720$ conics (mostly reducible), suggesting that $432$ should
be far from the maximum $\Nmax4(2)$. Therefore, in this note I suggest a very
simple (although also implicit) construction of a smooth quartic with $800$
irreducible conics.

\theorem[see \autoref{proof.main}]\label{th.main}
There exists a smooth quartic surface $\bX\subset\Cp3$ containing $800$
irreducible conics.
\endtheorem

{\3The quartic~$X_4$ is Kummer in the sense of~\cite{Barth.Bauer:conics,Bauer:conics}:
it contains $16$ disjoint conics.}
I conjecture that $\Nmax4(2)=800$ and, moreover,
$800$ is the sharp upper bound on the total number of
conics (irreducible or reducible) in a smooth spatial quartic.

{\3There has been a considerable development after this note appeared in the
\texttt{arXiv}.
X.~Roulleau observed that,
computing the projective automorphism group and using~\cite{Bonnafe.Sarti},
%we conclude
%%Xavier Roulleau concluded
%that the quartic
$\bX$ in \autoref{th.main}
must be given by the Mukai polynomial
\[*
z_0^4+z_1^4+z_2^4+z_3^4+12z_0z_1z_2z_3=0;
\]
even though only $320$ conics were found in~\cite{Bonnafe.Sarti}.
Then, B.~Naskr\k{e}cki found explicit equations of all $800$ conics.}

%\remark
%{\3There has been a considerable development after this note appeared in the
%\texttt{arXiv}.
%X.~Roulleau observed that,
%computing the projective automorphism group and using~\cite{Bonnafe.Sarti},
%%we conclude
%%%Xavier Roulleau concluded
%%that the quartic
%$\bX$ in \autoref{th.main}
%must be given by the Mukai polynomial
%\[*
%z_0^4+z_1^4+z_2^4+z_3^4+12z_0z_1z_2z_3=0;
%\]
%%even though only $320$ conics were found in~\cite{Bonnafe.Sarti}.
%%(this fact was discovered by X.~Roulleau),
%then, B. Naskrecki found all $800$ conics.
%(Only $320$ conics are presented~\cite{Bonnafe.Sarti}.)}
%\endremark

\subsection*{Acknowledgements}
I am grateful to S{\l}awomir Rams who introduced me to the subject and made
me familiar with~\cite{Barth.Bauer:conics,Bauer:conics}.
%the current state of the art.
%;
%it is his curiosity that encouraged my
%work in this direction.
{\3I extend my gratitude to Dino Festi,
Bartosz Naskr\k{e}cki,
and Xavier Roulleau,
who took active part in the
study of~$X_4$.}
%and drew my attention to~\cite{Bonnafe.Sarti}.

\section{The Leech lattice\pdfstr{}{ \rm(see \cite{Conway.Sloane})}}\label{S.Leech}

\subsection{The Golay code}\label{s.Golay}
The \emph{\rom(extended binary\rom) Golay code} is the only binary code of
length~$24$, dimension~$12$, and minimal Hamming distance~$8$. We regard
codewords as subsets of $\Omega:=\{1,\ldots,24\}$ and denote this collection
of subsets by~$\GC$; clearly, $\ls|\GC|=2^{12}$. The code~$\GC$ is invariant
under the complement $o\mapsto\Omega\sminus o$. Apart from~$\varnothing$
and~$\Omega$ itself, it consists of $759$ \emph{octads} (codewords of
length~$8$), $759$ complements thereof, and $2576$ \emph{dodecads} (codewords
of length~$12$).

The setwise stabilizer of~$\GC$ in the full symmetric
group $\Bbb{S}(\Omega)$ is the
Mathieu group $M_{24}$ of order $244823040$; the actions of this group
on~$\Omega$ and~$\GC$ are described in detail in \S\,2 of
\cite[Chapter~10]{Conway.Sloane}.

\subsection{The square~$4$ vectors}\label{s.Leech}
The \emph{Leech lattice} is the only root-free unimodular even
positive definite lattice of rank~$24$.
For the construction, consider the standard Euclidean
lattice $E:=\bigoplus_i\Z e_i$, $i\in\Omega$, and divide the form by $8$,
so that $e_i^2=1/8$. (Thus, we avoid the factor $8^{-1/2}$ appearing
throughout in~\cite{Conway.Sloane}.) Then, the Leech lattice is the
sublattice $\Lambda\subset E$ spanned over~$\Z$ by
the square~$4$ vectors of the form
\[
(\mp3,\pm1^{23})\qquad
\text{(the upper signs are taken on a codeword $o\in\GC$)}.
\label{eq.31}
\]
(We use the notation of~\cite{Conway.Sloane}: $a^m,b^n,\ldots$
means that there are $m$ coordinates equal to~$a$, $n$ coordinates equal
to~$b$, \etc.)
Apart from~\eqref{eq.31}, the square~$4$ vectors in~$\Lambda$ are
\begin{align}\label{eq.20}
(\pm2^8,0^{16})&\qquad
 \text{($\pm2$ are taken on an octad, the number of $+$ is even)},\rlap{ or}\\
(\pm4^2 ,0^{22})&\qquad\label{eq.40}
 \text{(no further restrictions)}.
\end{align}
Altogether, there are $196560$ square~$4$ vectors:
$24\cdot\ls|\GC|=98304$ vectors as in~\eqref{eq.31},
$2^7\cdot759=97152$ vectors as in~\eqref{eq.20}, and
$2^2\cdot C(24,2) = 1104$ vectors as in~\eqref{eq.40}.

\section{The construction}\label{S.proof}

In this section, we prove \autoref{th.main}.

\subsection{The lattice~$S$}\label{s.S}
Consider the lattice $\tV:=\Z\hbar+\Z a+\Z u_1+\Z u_2+\Z u_3$ with the Gram matrix
\[*
\def\-{\llap{$-$}}
\begin{bmatrix}
  4&  2&  0&  0&  0 \cr
  2&  4&  2&  0&  1 \cr
  0&  2&  4&  2&\-1 \cr
  0&  0&  2&  4&  0 \cr
  0&  1&\-1&  0&  4
\end{bmatrix}.
\]
It can be shown that, up to $\OG(\Lambda)$, there is a unique primitive
isometric embedding $\tV\to\Lambda$; however, for our example,
we merely choose a particular model.
Fix an ordered quintuple $\fixed:=(1,2,3,4,5)\subset\Omega$ and choose one of
the four octads~$\octad$ such that $\octad\cap\fixed=\{1,2,4,5\}$
(\cf. \emph{sextets} in \S\,2.5 of \cite[Chapter 10]{Conway.Sloane});
upon reordering~$\Omega$,
we can assume that $\octad=\{1,2,4,5,6,7,8,9\}$ (the underlined positions in
the top row of \autoref{tab.V}).
Then, the generators of~$\tV$ can be chosen as shown in the
upper part of \autoref{tab.V}. (For better readability, we represent zeros by
dots; all components beyond $\set:=\fixed\cup\octad$ are zeros.)

The choice of~$\fixed$ and~$\octad$ is unique up to~$M_{24}$; furthermore,
the subgroup $G\subset M_{24}$ stabilising~$\fixed$ pointwise and~$\octad$ as
a set can be identified with the alternating group
$\Bbb{A}(\octad\sminus\fixed)$; in particular, it acts simply transitively on
the set of ordered pairs
\[
(p,q)\:\quad p,q\in\octad\sminus\fixed=\{6,7,8,9\},\quad p\ne q.
\label{eq.pq}
\]

\table
\caption{The lattice~$\tV$ and the conics}\label{tab.V}
\hbox to\hsize{\hss
 \def\kk{\kern2ex\relax}%
 \def\bx#1{\kk\hss$#1$\hss\kk}%
 \def\vr{\vrule width0pt}
 \let\.\cdot\catcode`\.\active\let.\.
 \let\--\catcode`\-\active\def-{\llap{$\-$}}%
 \let\**\catcode`\*\active\def*{\rlap{$^\*$}}%
 \let\s\scriptstyle
 \def\u#1{\s\underline{#1}}
 \vbox{\offinterlineskip\ialign{\hss\,\,$#$\hss\quad\vrule\,&
  \bx{#}&\bx{#}&\bx{#}&\bx{#}&\bx{#}\vrule&
   \bx{#}&\bx{#}&\bx{#}&\bx{#}\vrule&\,\,\,$#$\,\strut\hss\cr
     \s\#&\u1&\u2&\s3&\u4&\u5&\u6&\u7&\u8&\u9\vr depth6pt height7pt\cr
    \hbar&  4&  4&  .&  .&  .&  .&  .&  .&  .&\cr
        a&  .&  4&  4&  .&  .&  .&  .&  .&  .&\cr
      u_1&  .&  .&  4&  4&  .&  .&  .&  .&  .&\cr
      u_2&  .&  .&  .&  4&  4&  .&  .&  .&  .&\cr
      u_3& -2&  2&  .& -2&  2&  2&  2&  2&  2&\vr depth4pt\cr
       \noalign{\hrule}
    \num1&  1&  3& -1&  1& -1&  1&  1&-1*&-1*&\pm1^{15}\vr height10pt\cr
    \num2&  3&  1&  1& -1&  1&  1&  1&-1*&-1*&\pm1^{15}\cr
    \num3&  2&  2&  .&  .&  .&  .&  .&  .&  .&\pm2^6,0^9\cr
    \num4&  2&  2&  .&  .&  .&  .&  .& 2*&-2*&\pm2^4,0^{11}\cr
         &\multispan5\hss\small $\text{fixed}=\fixed$\hss\vrule&
          \multispan4\hss\small movable in $\octad\sminus\fixed$\hss\vrule depth2pt height8pt\cr
\crcr}}\hss}
\endtable
\table
\caption{The number of conics in $S$}\label{tab.S}
\hrule height0pt\relax\vskip-\abovedisplayskip
\[*
\def\0{\phantom{0}}
\alignedat3
\num1\:&\ & C(4,2)\cdot\underline{16}&=\096&\quad&
 \text{(codewords $o\in\GC$ such that $o\cap\set=\{2,3,5,p,q\}$)},\\
\num2\:&& C(4,2)\cdot\underline{16}&=\096&&
 \text{(codewords $o\in\GC$ such that $o\cap\set=\{1,4,p,q\}$)},\\
\num3\:&& 2^5\cdot\underline{10}&=320&&
 \text{(octads $o\in\GC$ such that $o\cap\set=\{1,2\}$)},\\
\num4\:&& 2^3\cdot P(4,2)\cdot\underline{3}&=288&&
 \text{(octads $o\in\GC$ such that $o\cap\set=\{1,2,p,q\}$)}.
\endalignedat
\]
\vskip-\belowdisplayskip
\endtable

Define a \emph{conic} as a square~$4$ vector $l\in\Lambda$ such that
\[*
l\cdot\hbar=2,\quad l\cdot a=1,\quad l\cdot u_1=l\cdot u_2=l\cdot u_3=0.
\]
This strange condition can be recast as follows: $l\cdot\hbar=2$ and $l$
(as well as $\hbar$) lies in the rank~$20$ lattice
\[*
S:=\bV^\perp\subset\Lambda,\quad
\text{where $\bV:=\hbar^\perp_\tV$}.
\]
Using \autoref{s.Leech}, we conclude that each conic fits one of the four
patterns shown at the bottom of \autoref{tab.V}: there are two
for~\eqref{eq.31} and two for~\eqref{eq.20}. (If $l$ is as in~\eqref{eq.40},
we have $l\cdot a=0\bmod2$.)
The number of conics within each pattern is
computed as shown in \autoref{tab.S}, where
\roster*
\item
the ordered or unordered pair $(p,q)$ as in~\eqref{eq.pq} designates the two
variable special positions marked with a $^*$ in \autoref{tab.V},
\item
the underlined factor
counts certain codewords $o\in\GC$; the restrictions given by~\eqref{eq.31}
or~\eqref{eq.20} are described in the parentheses, and
\item
the other factors account for the choice of $(p,q)$ and/or signs in $\pm2$.
\endroster
These counts sum up to $800$.

\subsection{The N\'{e}ron--Severi lattice}\label{s.NS}
Observe that
$\hbar\in2S\dual$: indeed,
$\hbar-2a\in\bV$ and we have $x\cdot\hbar=2x\cdot a=0\bmod2$ for any
$x\in S$.
Thus, we can apply to $S\ni\hbar$ the construction of~\cite{degt:conics},
\ie, consider the orthogonal complement
$\hbar^\perp_S=\tV^\perp\subset\Lambda$, reverse the sign of the form,
and pass to the index~$2$ extension
\[*
N:=\bigl(-(\hbar^\perp_S)\oplus\Z h\bigr)^\sim_2,\quad h^2=4,
\]
containing the vector $c:=c(l):=l-\frac12\hbar+\frac12h$ for some
(equivalently, any) conic $l\in S$.
These $800$ new vectors $c\in N$ are also called
\emph{conics}; one obviously has
\[
c^2=-2\quad\text{and}\quad c\cdot h=2.
\label{eq.conic}
\]
%They are in a bijection with the conics in~$S$; hence, there are $800$ of
%them.

%\remark
%In order to use~\cite{degt:conics}, we need that
%$\hbar\in2S\dual$. This follows from the fact that
%$\hbar-2a\in \bV$ and, hence, $x\cdot\hbar=2x\cdot a=0\bmod2$ for any
%$x\in S$.
%\endremark

Starting from
\[*
\discr \tV\cong\bmatrix
1&\frac12\cr
\frac12&1
\endbmatrix
 \oplus\bmatrix\frac18\endbmatrix
 \oplus\bmatrix\frac25\endbmatrix
\]
(see Nikulin~\cite{Nikulin:forms} for the concept of \emph{discriminant form}
$\discr \tV:=\tV\dual\!/\tV$ and related techniques),
we easily compute
\[*
\Cal{N}:=\discr N\cong\bmatrix\frac54\endbmatrix
 \oplus\bmatrix\frac18\endbmatrix
 \oplus\bmatrix\frac25\endbmatrix
\cong\bmatrix-\frac14\endbmatrix
 \oplus\bmatrix-\frac58\endbmatrix
 \oplus\bmatrix\frac25\endbmatrix.
\]
Therefore, $-\Cal{N}\cong\discr T$, where
$T:=\Z b\oplus\Z c$, $b^2=4$, $c^2=40$.
Then, it follows from~\cite{Nikulin:forms} that there is a primitive
isometric embedding of the hyperbolic lattice $N$
to the intersection lattice $H_2$ of a
$K3$-surface,
so that $T\cong N^\perp$ plays the r\^{o}le of the transcendental lattice.
Finally, by the surjectivity of the period
map~\cite{Kulikov:periods},
we conclude that there exists a $K3$-surface~$\tX$ with
$\NS(\tX)\cong N$.

\subsection{Proof of \autoref{th.main}}\label{proof.main}
The N\'{e}ron--Severi lattice $\NS(\tX)\cong N$ constructed in the previous
section is equipped with
a distinguished polarisation $h\in N$, $h^2=4$. Since the
original lattice $S\subset\Lambda$ is root free, $N$ does \emph{not} contain
any of the following ``bad'' vectors:
\roster*
\item
$e\in N$ such that $e^2=-2$ and $e\cdot h=0$ (\emph{exceptional divisors}) or
\item
$e\in N$ such that $e^2=0$ and $e\cdot h=2$ (\emph{$2$-isotropic vectors})
\endroster
(see \cite{degt:conics} for details). Hence, by Nikulin~\cite{Nikulin:Weil}
and Saint-Donat~\cite{Saint-Donat}, the linear system $\ls|h|$ is fixed point
free and maps~$\tX$ onto a smooth quartic surface $\bX\subset\Cp3$.

The lattice $N$ contains $800$ conics~$c$ as in~\eqref{eq.conic}. By the
Riemann--Roch theorem, each class~$c$ is effective, \ie, represented by
a curve $C\subset\bX$ of projective degree~$2$. Since~$X$ is smooth and
contains no lines (or other curves of odd degree, as we have $h\in2N\dual$
by the construction), each of these curves~$C$ is irreducible.
This concludes the proof of \autoref{th.main}.
\qed

{
\let\.\DOTaccent
\def\cprime{$'$}
\bibliographystyle{amsplain}
\bibliography{degt}
}

\end{document}